**Werner DePauli-Schimanovich**
**Institute for Information Science, Dept. DB&AI, TU-Vienna, Austria**
**Werner.Depauli@gmail.com**

# A Brief History of Future Set Theory, K50-Set4

$$\forall \text{wff } B_j^0(x) \subseteq A_i^0(x): \text{Set}(\{|x: B_j(x)|\}) \not\vdash_{CT} \wedge \longrightarrow \text{Set}(\{|x: A_i(x)|\}).$$

(Werner DePauli, 1990)

**(0) Abstract[1]**
Mathematicians still use Naive Set Theory when generating sets without danger of producing any contradiction. Therefore their working method can be considered as a consistent inference system with an experience of over 100 years. My conjecture is that this method works well because mathematicians use only those predicates to form sets, which yield closed hereditary consistent predicate extensions. And for every open formula they use in the process of constructing of a certain (special) set (bottom up), we can always find an "almost-closed" formula (i.e. a parameter-free formula with only the free variable "x") which yields the same certain (special) set as predicate extension as constructed in the bottom up process before. Therefore the use of predicates with free parameters in the Comprehension Scheme does not cause any difficulties and can be "lifted" by meta-mathematical considerations.

KEYWORDS: naive set theory, Quine, new foundation, NF, universal sets, comprehension schema, predicate extension, philosophy of set theory, Zermelo, Fraenkel, ZF, complement.
CT is the Class-Theoretical frame: logic + "=" + Church Schema + Extensionality Axiom.

**(1) Introduction**
When David Hilbert delivered his famous talk at the 1928 International Congress of Mathematicians in Bologna, he claimed: "Cantor has created a Paradise for us from which nobody can expel us again!" This paradise, the Naïve "Mengenlehre", was based mainly on a single principle, which was taken as a dogma: "Sets are Predicate-Extensions." Or in other words: "For every arbitrary predicate, its extension (that is, all objects with this property) forms a set." After formalization this principle is called the general Comprehension Scheme (CoS) := forall wff A: forall y: [y in {x: A(x)} <==> A(y)].
"wff" is the abbreviation for well-formed formula.

---

[1] I want to thank Martin Goldstern, Randall Holmes, Matthias Baaz and Thomas Forster for their help and support by writing this paper.

In 1903, Bertrand Russell discovered the paradox[2] named after him, to wit, that the predicate "x non-in x" is a substitution-instance for A which falsifies the comprehension scheme. He wrote a letter to Gottlob Frege and told him his discovery. Frege had just finished his book "Grundlagen der Arithmetik" in which he wanted to put arithmetic and "Mengenlehre" onto a safe and sound fundation. He was forced to add an appendix: "Nothing worse can happen to a scientific writer than when his fundation breaks away after finishing his house."

Later other inconsistencies were found: circa 1908 Burali-Forti's paradox of the "set of all ordinals", circa 1909 Mirimanoff's paradox of the "set of all sets not containing an infinite descending element sequence", etc. In all these cases one single formula A was enough to falsify (CoS). Therefore set theoreticians posed themselves the question whether one should not restrict set formation to the "consistent" predicate-extensions (with free parameters); i.e. only accept those (non-pathological) A which, when substituted in (CoS), do not produce a contradiction. But soon it turned out that neither the complement "ko(m)"[3] of an arbitrary set "m" nor the Anti-Russell "x in x" allows one to deduce a contradiction if they are substituted into (CoS); but assuming both together simultaneously does so.

Therefore experts in set theory asked themselves whether the closed consistent predicate-extensions in (CoS) would yield a contradiction-free axiomatization of Naïve "Mengenlehre". (With "closed" formulas we mean in the following context always "almost-closed", i.e. formulas with the only free variable "x", over which we compress the predicate to a set. By this convention our predicate-extensions are always closed.) But this restriction was not enough. Because if you take some arbitrary undecideable formula[4] UD together with the implication of the Russell-condition "x non-in x", then its negation non-UD together with the implication of Russell is also undecideable. Each of these 2 implications substituted for A in (CoS) gives a consistent predicate-extension, but either UD or non-UD must be valid and this therefore produces the Russell-class as a set; a contradiction once again.

Now the ultimate question arises: do the "closed hereditary consistent predicate-extensions", in short CHC-PE (together only with the common axiom of extensionality (EE) and no other specific axioms) produce a contradiction-free collection of formulas?[5] I conjecture: "yes!", and until now nobody could show the contrary. I was also informed by Randell Holmes that he supports my opinion to demand hereditary-consistent[6] formulas A (instead of decidable ones). If this system of set theory is consistent, it is probably some version of positive set theory. This will give rise to a series of interesting investigations. Therefore we redefine "consistent" as "hereditary-consistent".[7]

---

[2] In the paper K50-Set3b in this book I analyzed the concept "paradox" in which I show that it can be split up into the 3 different notions "pathological", "antinomical" and "abnormal", each of which can be characterized formally.

[3] In the following we will always use K as abbreviation for complement of classes [as the German word "Komplement"], because the letter C is already used for Comprehension, Church, Cardinal, Conditional, Cofinal, etc, and k if the complement yields a set.

[4] A closed wff A shall be decidable in respect to FOL (with the only predicate "in") and (EE). If A is almost-closed (with only free x) we call it decidable if all its substitution instances with closed terms are decidable. Any PE shall be called decidable if its substituted A is it.

[5] This collection needs not to be recursively enumerable and therefore no (classical) axiom-system.

### (2) Separation Principle and Limitation of Size

At the time when this ultimate question arose (whether all CHC-PE's together establish a consistent domain of formulas), its investigation was already too late, because the fact that such predicate extensions have to be closed has been denied. Nobody wanted parameter-free set-existence schemata. Georg Cantor had discovered in 1899 the antinomy named after him, i.e. that on the one hand the universal set "us" is of smaller power then its power-set "p(us)": [Card(us) < Card(p(us))], but on the other hand the power-set would have to be extensionally equal to the universal set: [us = p(us)]. This is a classical contradiction, because there exists of course a function from "us" onto "us" (i.e. the identity-function).

This antinomy follows directly from the theorem of Cantor, which says that every set is of smaller power than its power-set. But the theorem of Cantor can be applied to the universal set only because the „Aussonderungs-Prinzip" (= principle of separation) for Cantor was something like an "allgemeingültiges" (generally valid) axiom for arbitrary sets! Therefore it would not have made any sense to consider the closed consistent predicate extensions, because the set theoreticians needed parameters like those used in the separation-principle.

Today the most common separation principle is the axiom-scheme (Sep):
forall wff A: forall x: Set({y: y in x & A(y)}). Or in another formulation:
exist y: forall x: [y = (x intersection A)], where A is the class {|z : A (z)|} of all objects z with A.[8] If we want to accept the principle of separation as a generally valid axiom for arbitrary sets (which can occur as free parameters), then the „Limitation of Size" ideology[9] follows directly. Therefore Hausdorff and likewise Fraenkel wrote in their fundamental textbooks (in the 1920s): "The antinomies arise through the construction of sets which are extensionally equal with the universe!" Based on this ideology is the axiomatic system of set theory (today most well-known), Zermelo-Fraenkel (= ZF), and adding the axiom of choice (AC), the system ZFC.

---

[6] Hereditary-Consistent (= HC) means: Restrict comprehension to the formulas with the property that all subformulas (suitably defined) are consistent (i.e. non-pathological). This yield the CHC-PE.

[7] Concerning set-existence the CHC-PE's are a "maximal consistent set of formulas" (in a certain sense), and the corresponding sets form a "maximal model".

[8] In class theory (Sep) is the subset axiom: Every subclass of a set is again a set: forall X: [X =< y ==> Set (X)].

[9] See Hallet [1984]: Cantorian Set Theory and the Limitation of Size.

592

Today this "Limitation of Size" ideology cannot be kept any longer. In 1937 Willard Van Orman Quine published his system "New Foundation" (= NF)[10], which he developed from Russell's theory of types. This system contains a universal set and must thus be inconsistent in virtue of Fraenkel. But in 1969 Ronald Björn Jensen[11] showed the relative consistency of NF with ZFU (= ZF with Urelemente). Since the days of Zermelo-Fraenkel a lot of set theoretical systems with a universal set have been developed in the meantime. A beautiful collection of these systems can be found in the book of Thomas Forster.[12]

**(3) Small Sets, Complementation, and ZF**
In my dissertation, "Extension der Mengenlehre" 1971 (K50-Set1 in this book), and its popular-scientific edition "Der Mengenbildungs-Prozess" (1971 too, in the journal for literature „Manuskripte '33" [or K50-Set2 in this book]) I had raised the question as to a consistent axiomatization of the Naïve Mengenlehre. At that time I could present only partial results, e.g. the answer to the question "What is the slightest modification of ZF needed to add a complement-axiom?" My answer was: the restriction of the axiom of separation to "small" sets. Of course, the axiom of replacement -- which says, formulated in the language of class theory, that "Every image of a set is again a set" -- has to be restricted too. (In fact it's enough to restrict replacement, because separation is a logical consequence of it.)

The vague notion of a "small" set can be made precise in several ways:
(1) Small shall be represented by well-foundedness, or hereditary foundedness,
(2) The Cantorian sets (for which the theorem of Cantor is valid) are the small ones,
(3) "Slim" or "Komplements-schmächtiger", i.e. the small sets are those with smaller power than their complement.

In the simple language of class theory, the restricted axiom of separation says: „The subclass of every small set is again a set." The analogue is valid for the images. No mathematician would ever try to construct a subset or image of the universal set, or of some other set "extensionally equal to the universe". In any case, the mathematicians always do the right thing!

If we restrict the axiom of replacement to small sets (and from this follows already the restricted axiom of separation), then this is (at least for me) the slightest modification of ZF which allows us to add the existence of a general complement. In 1971 at the International Logic-Congress in Bukarest I gave a talk on this subject and explained it to Dana Scott (and after that I discussed it with him in Vienna). At that time I even hoped to obtain a proof of the relative consistency of my system with ZF. But Martin Goldstern convinced me 5 years ago that this proof is false. (And good-minded as I am, I did believe him.) He published another proof for a stronger weakened ZF on his homepage.[13]

---

[10] A good introduction to NF the reader can find in Holmes [1998] in the References.
[11] See Jensen [1969] in the References K50-Set10.
[12] See Forster [1995] in the References K50-Set10.
[13] See Goldstern [1998] in the References K50-Set10.

After all which has been said until now, one fact should have become evident: If there exists a universal set or the complement in a set theoretical system, then the axioms of replacement and separation have to be restricted to small sets. Cantor wanted the general separation, but that guided us directly into the dead-end street of the Limitation of Size Ideology and therefore to ZF. Despite the fact that ZF is the favorite system of logicians, nevertheless, it is highly unnatural. No mathematician constructs his sets bottom-up by iterated application of the ZF-axioms!

But this was not the reason why ZF was created. Concerning Hilbert's Program the scientific community wanted to show the consistency of set theory. But since Kurt Gödel we know that this is impossible, because (for a sufficiently rich system) one cannot show its consistency by means of the system itself. Therefore ZF has lost its original justification, and it is today in common use mainly because it has been spread all over the world and it is somehow easy to handle (compared e.g. with NF). Mathematicians (especially Bourbaki[14]) have in any case always ignored ZF, and they are still working with the Naïve "Mengenlehre", which they use contentedly without generating any contradictions.

**(4) How do Mathematicians Generate Sets?**
Taking the logical point of view which I have described just above (i.e. if somebody accepts the restriction of replacement and separation), naturally the question arises how the working method of the mathematicians can be justified philosophically? (ZF certainly does not do the job!) Therefore I shall investigate here the methods mathematicians actually use when they are generating sets.

Working mathematicians have a very strong intuition. And this intuition gives them the insight that a formula A does not produce a contradiction in the comprehension scheme (CoS) even when it is very complicated. They also use only hereditary-consistent formulas and if possible no undecidable ones. (Since these formulas are usualy very large [like the Gödel-formula in Gödel's proof or the theorem of Paris-Harrington] mathematicians never use them in practise). Of course, mathematicians often use open decidable and hereditary-consistent predicate extensions with free parameters, e.g. for constructing the union or intersection of finitely many sets, or n-tupels. But this problem can be lifted once we agree that the CHC-PE's (together with extensionality axiom (EE), that two sets containing the same elements are equal) can be considered to be known as contradiction-free by experience. Then the use of such open extensions can be shown to be correct by meta-theorems stating that in such cases, a corresponding CHC-PE can also be constructed which forms the same set.[15]

Mathematicians only make relative complements (e.g., in a sigma algebra). Therefore they will not feel disturbed if it were allowed to construct absolute complements of small sets, too. This can already be deduced from the CHC-PE's yet. Also the axiom of choice is applied by the mathematicians only to small sets. But a restricted (to small sets) axiom of choice (SmallAC) is also compatible with Quine's NF, and this gives rise to (and hope for) the assumption that it will also be consistent together with the CHC-PE plus (EE).

---

[14] See Mathias [1992] "The Ignorance of Bourbaki" in K50-Set10.
[15] A great help in this case would be if we add the term-consistency rule: If for every closed term t of a wff A the statement A(t) holds, then "forall x: A(x)".

Therefore if mathematicians take a set-constituting property A and substitute it into the set operator, it always produces a set, and no contradiction. Because mathematicians only use hereditary-consistent formulas A (which make the PE's also as a collection consistent), no contradiction in practise arises. If they use formulas A, with only "x" as free variable, substituting them into the set operator, it makes our philosophical justification easier. But if they use also formulas with free parameters in turn, it does not create difficulties as long as they use only small sets (where the pathological ones are excluded automatically). This can be considered as a 100-year-long experience of consistency of naïve set theory in practise![16] (That is several years more than ZF or NF exists.)

### (5) The Systems NAM* and NACT*

The collection of all closed hereditary-consistent predicate extensions (= CHC-PE) together with (EE) form a new system of set theory. But since this collection of formulas is not recursively enumerable, it is not a axiom system in the classical sense. Anyway, the property of a formula A to establish a CHC-PE is easy to grasp, and therefore it does not matter whether the set of all CHC-PE's is not recursively enumerable, because it is enough to consider some large recursively enumerable subset of it. Let us therefore call this system of all CHC-PE's plus (EE) an axiom-system in the wider sense; we want to name it NAM* (= Naïve Axiomatic Mengenlehre Star), or if we formalize it in the more simple language of class theory, NACT* (= Naïve Axiomatic Class Theory Star, or: the naked star).

NAM* works in the following way: First we want to find the wffs A with only one free variable x which itself (and all its subsets) do not produce a contradiction substituted into (CoS). Let us therefore consider a formula generator FG which enumerates all wffs $A_j$ in a natural way. For all $A_j$ call their substitution instances (CoS-$A_j$). A decision algorithm DA shall find out if (CoS-$A_j$) together with (EE) produces a contradiction only by use of the rules of pure logic (with the only predicate "is an element of") and without any special axioms of set theory at all! If this happens we want to call such a predicate $A_j$ "pathological".[17] Therefore the hereditary-non-pathological $A_j$'s are exactly the CHC-PE's. NAM* is therefore (CoS) restricted to the (hereditary-)non-pathological $A_j$'s, plus (EE).
As a formula: NAM* := (CoS) & (EE), where (CoS)={(CoS-$A_j$)} $j \epsilon$ J with HnP($A_j$).

---

[16] There must exist an adequate formalization of the contentually working method of the mathematicians using naïve set theory. The question is only: Who will find it first? Of course: implicitely I conjecture in this paper: The CHC-PE's are this adequate formalization of mathematical practise.

[17] Of course this DA proves also if non-(CoS-$A_j$) produces a contradiction, and therefore if $A_j$ [or better with some formula B its corresponding PE] is satisfyable or not.

For our purpose (as already mentioned) it does not matter that such a decision algorithm DA does not exist for all (hereditary-)non-pathological predicates. It's enough to find one for a sufficiently large class of non-Patho, which should signify the (second order class of all) pathological classes. The existence of such a DA for all (CoS-Aj), with the intention of producing the whole class Consi [= Hereditary non-Patho], is logically equivalent to the decision problem of FOL, shown to be impossible to solve by Alonzo Church 1936. The analogue statement is valid also for NACT* (= the naked star): classes with (hereditary) non-patho wffs (i.e. the CHC-PE's) are sets.

Therefore we can conclude that NAM* (i.e. the most important system instantiating NAM) is a real naïve set theory, very similar to the one used by mathematicians in practice, with a unique logical principle to generate sets (on which it is based), in contrast to other systems with more or less arbitrarily selected axioms. This shows also that Logicism is not dead yet. NAM* is a rectification of mathematical practise, and the systems ZF and NF are (in a certain meta-mathematical sense) included in NAM*. (Of course we have to consider in this case the restricted axioms of separation, replacement and choice, which have to be shown as meta-theorems, and consider only the parameter-free statements.)

Set theory was already nearly dead, after all the ZF people had done with it. Dana Scott said on the occasion of receiving the Bolzano award in Prague: "Bernhard Bolzano would have been very disappointed if he had known what excesses set theory of today is producing!"
    Now, after creating NAM* and NACT*, it will become interesting again. Because Naïve Axiomatic "Mengenlehre" in general (which we want to call short hand NAM in the following) delivers in addition an instrument to experiment with set theory. This way we can find out what mathematicians consider to be normal sets. This is a philosophical program like the one in the '30s of the past century to formalize the intuitive notion of computability. NAM* is only the most significant system found by this experimentation. In the article "Naïve Axiomatic Mengenlehre for Experiments" we want to show how experimentation with NAM works in the following. (See K50-Set7b).

**References for the Set-Theoretical Articles, K50-Set10**
Bernays, Paul & Fraenkel, Abraham
[1968] Axiomatic Set Theory. Springer Verlag, Heidelberg, New York.
Brunner, Norbert & Felgner, Ulrich
[2002] Gödels Universum der konstruktiblen Mengen. In: [Buldt & al, 2002].
Buldt, Bernd & Köhler, Eckehard & Stöltzner, Michael & Weibel, Peter & Klein, Carsten & DePauli-Schimanovich-Göttig, Werner
[2002] Kurt Gödel: Wahrheit und Beweisbarkeit, Band 2: Kompendium zu Gödels Werk. oebv&htp, Wien.
Casti, John & DePauli, Werner
[2000] Gödel: A Life of Logic. Perseus Publishing, Cambridge (MA).
Davidson, Donald & Hintikka, Jaakko
[1969] Words and Objections. D. Reidel Publ. Comp., Dordrecht.

DePauli-Schimanovich, Werner (See also: Schimanovich [1971a], [1971b]) and [1981]).
[ca 1990] On Frege's True Way Out. Article K50-Set8b in this book.
[1998] Hegel und die Mengenlehre. Preprint at: http://www.univie.ac.at/bvi/europolis .
[2002a] Naïve Axiomatic Mengenlehre for Experiments. In: Proceedings of the HPLMC (= History and Philosophy of Logic, Mathematics and Computing), Sept. 2002 in Hagenberg/Linz (in honour of Bruno Buchberger).
[2002b] The Notion "Pathology" in Set Theora. In: Abstracts of the Conference HiPhiLoMaC (= History and Philosophy of Logic, Mathematics and Computing), Nov. 2002 in San Sebastian.
[2005a] Kurt Gödel und die Mathematische Logik (EUROPOLIS5). Trauner Verlag, Linz/Austria.
[2005b] Arrow's Paradox ist partial-consistent. In: DePauli-Schimanovich [2005a].
[2006a] A Brief History of Future Set Theory. Article K50-Set4 in this book.
[2006b] The Notion of "Pathology" in Set Theory. Article K50-Set5b in this book.
[2006c] Naïve Axiomatic Class Theory NACT: a Solution for the Antinomies of Naïve "Mengenlehre". Article K50-Set6 in this book.
[2006d] Naïve Axiomatic Mengenlehre for Experiments. Article K50-Set7b in this book.
DePauli-Schimanovich, Werner & Weibel, Peter
[1997] Kurt Gödel: Ein Mathematischer Mythos. Hölder-Pichler-Tempsky Verlag, Wien.
Feferman, Solomon & Dawson, John & Kleene, Stephen & Moore, Gregory & Solovay, Robert & Heijenoort, Jean van,
[1990] Kurt Gödel: Collected Works, Volume II. Oxford University Press, New York & Oxford.
Felgner, Ulrich
[1985] Mengenlehre: Wege Mathematischer Grundlagenforschung. Wissensch. Buchgesellschaft, Darmstadt.
[2002] Zur Geschichte des Mengenbegriffs. In: [Buldt & al, 2002] und Artikel K50-Set3 in this book.
Forster, Thomas
[1995] Set Theory with an Universal Set. Exploring an Untyped Universe. (2$^{nd}$ Edition.) Oxford Science Publ., Clarendon Press, Oxford.
Gödel, Kurt
[1938] The relative consistency of the axiom of choice and of the generalized continuum hypothesis. In: [Feferman & al, 1990].
Goldstern, Martin & Judah, Haim
[1995] The Incompleteness Phenomenon. A. K. Peters Ltd., Wellesley (MA).
Goldstern, Martin
[1998] Set Theory with Complements. http://info.tuwien.ac.at/goldstern/papers/notes/zfpk.pdf
Hallet, Michael
[1984] Cantorian Set Theory and Limitation of Size. Oxford Univ. Press, N.Y. & Oxford.
Halmos, Paul
[1960] Naive Set Theory.   Van Nostrand Company Inc., Princeton (NJ).
Holmes, Randall
[1998] Elementary Set Theory with a Universal Set.

Vol. 10 of the Cahiers du Centre de Logique. Academia-Bruylant, Louvain-la-Neuve (Belgium).
[2002] The inconsistency of double-extension set theory.
http://math.boisestate.edu/~holmes/holmes/doubleextension.ps
Jech, Thomas
[1974] Procedings of the Symposium in Pure Mathematics (1970), Vol. XIII, Part 2, AMS, Provicene R.I. .
Jensen, Ronald Björn
[1969] On the consistency of a slight (?) modification of Quine's New Foundation. In: [Davidson & Hintikka, 1969].
Köhler, Eckehart & Weibel, Peter & Stöltzner, Michael & Buldt, Bernd & Klein, Carsten & DePauli-Schimanovich-Göttig, Werner
[2002] Kurt Gödel: Wahrheit und Beweisbarkeit, Band 1: Dokumente und historische Analysen. Hölder-Pichler-Tempsky Verlag, Wien.
Kolleritsch, Alfred & Waldorf, Günter
[1971] Manuskripte 33/'71 (Zeitschrift für Literatur und Kunst). Forum Stadtpark, A-8010 Graz, Austria.
Mathias, A.R.D.
[1992] The Ignorance of Bourbaki. In: The Mathematical Intelligencer 14 (No.3).
Quine, Willard Van Orman
[1969] Set Theory and its Logic. Belknap Press of Harvard University Press, Cambridge (MA).
Rubin, Jean & Rubin, Herman
[1978] Equivalents of the Axiom of Choice.    Springer, Heidelberg & New York
Schimanovich, Werner,
[1971a] Extension der Mengenlehre. Dissertation an der Universität Wien.
[1971b] Der Mengenbildungs-Prozess. In: [Kolleritsch & Waldorf, 1971].
[1981] The Formal Explication of the Concept of Antinomy. In: EUROPOLIS5, K44-Lo2a (and Lo2b). Or: Wittgenstein and his Impact on Contemporary Thought, Proceedings of the 2$^{nd}$ International Wittgenstein Symposium (29$^{th}$ of Aug. to 4$^{th}$ of Sept. 1977), Hölder-Pichler-Tempsky, Wien 1978.
Schimanovich-Galidescu, Maria-Elena
[2002] Princeton – Wien, 1946 – 1966. Gödels Briefe an seine Mutter. In: [Köhler & al, 2002].
Scott, Dana,
[1974] Axiomatizing Set Theory. In: Jech [1974] .
Suppes, Patrick,
[1960] Axiomatic Set Theory. D. Van Nostrand Company Inc., Princeton (NJ).
Weibel, Peter & Schimanovich, Werner
[1986] Kurt Gödel: Ein mathematischer Mythos. Film, 80 minutes, copyright ORF (= Austrian Television Network), ORF-shop, Würzburggasse 30, A-1130 Wien.
Wittgenstein, Ludwig,
[1956] Bemerkungen über die Grundlagen der Mathematik / Remarks on the Foundation of Mathematics. The M.I.T. Press, Cambridge MA and London.